\documentclass[reqno]{amsart}		
\usepackage{amsmath}
\usepackage{amssymb}
\usepackage{eucal}
\usepackage{enumerate}

\renewcommand\epsilon{\varepsilon}
\renewcommand\phi{\varphi}

\newcommand\LieAlg[1]{\mathfrak{#1}}
\newcommand\field[1]{\mathbb{#1}}
\newcommand\category[1]{\mathfrak{#1}}

\newcommand\R{{\field R}}
\newcommand\C{{\field C}}

\newcommand\set[1]{\left\{ #1 \right\}}
\newcommand\setst[2]%
	{\left\{\,#1\vphantom{#2} \;\right|\left. #2 \vphantom{#1}\,\right\}}

\newcommand\Ad{\operatorname{Ad}\nolimits}
\newcommand\ad{\operatorname{ad}\nolimits}

\newcommand\Id{\operatorname{Id}}
\renewcommand\Im{\operatorname{Im}}

\newcommand\from{\colon}

\newcommand\inv{^{-1}}
\newcommand\comp{\mathbin{\raise.29ex\hbox{$\scriptscriptstyle\circ$}}}
\newcommand\contr{\mathbin{\raise.29ex\hbox{$\lrcorner$}}}
\newcommand\blank{\makebox[0.8em]{$\cdot$}}

\newcommand\iso{\cong}

\newcommand\tensor{\otimes}

\newcommand\defeq{\overset{\text{\rm def}}{=}}
\newcommand\intersect{\cap}

\newcommand\G{\category G}
\newcommand\g{\LieAlg g}
\newcommand\h{\LieAlg h}
\newcommand\p{\LieAlg p}

\newcommand\Xitwiddle{\tilde\Xi}
\newcommand\thetabar{\bar \theta}
\renewcommand\L{{\mathcal L}}
\renewcommand\i{\iota}
\renewcommand\O{{\mathcal O}}	

\newtheorem{theorem}{Theorem}
\newtheorem{proposition}{Proposition}
\theoremstyle{definition}
\newtheorem*{definition}{Definition}
\theoremstyle{remark}
\newtheorem*{remark}{Remark}
\newtheorem*{remarks}{Remarks}
\newtheorem{lemma}{Lemma}

{\begin{small}[\relax}{]\end{small}}

\title{Symmetric Pairs and Moment Spaces}
\author{Matthew Leingang}
\date\today
\address{Harvard University Department of Mathematics, Cambridge, MA
02138-2109, U.S.A.}
\email{leingang@math.harvard.edu}
\begin{document}
\begin{abstract}
For a Lie group $G$, we seek the right definition of a {\em moment
space} for $G$.  One axiom is clear, involving a closed equivariant
three-form.  We construct this form for symmetric spaces associated to
symmetric pairs $(H,G)$ with an additional structure.  Furthermore, we
prove a decomposition theorem for these pairs over a compact, connected,
and semisimple group.

This monograph details work in progess and is not complete.
\end{abstract}
\maketitle

\addtocounter{section}{-1}
\section{Introduction}

The concept of momentum in a physical system with symmetry is quite
classical.  The collection of conserved quantities becomes in a more
modern viewpoint a map from phase space to the dual of the Lie algebra
of the symmetry group.  When such a {\em moment map}\/ exists, much
more can be said about the topology of phase space.  This is the
well-known theory of Hamiltonian $G$-spaces.

Recently, though, Alekseev, Meinrenken, and Malkin \cite{AMM} have
formalized a theory of ``$q$-Hamiltonian'' $G$-spaces, in which the
momentum is $G$- rather than $\g^*$-valued.  There is even a third
example for target of a moment map, which corresponds to a noncompact
symmetric space.

This leads us naturally to the question of the number of spaces which
can serve as target of a moment map.  We provide a recipe which given
a symmetric space with structure group $G$ and a special pairing on
the Lie algebra of the associated symmetric pair, produces a moment
space for $G$.  We also show that for suitable $G$ the only moment
spaces which arise in this manner are the known ones and combinations
of them.  This includes results proven independently by Alekseev,
Malkin, Meinrenken, and Kosmann-Schwarzbach (\cite{A+K-S} and
\cite{AMM-MP}).  

What remains to be shown is that every moment space can be produced in
this manner.  We hope that this the addition of further conditions to
the definition of moment space (corresponding to other aspects of the
classical Hamiltonian and $q$-Hamiltonian momenta), this can be done.

Throughout we will assume $G$ is compact.  By ``$G$-space'' we will
always mean manifold or sometimes orbifold on which $G$ acts.  Given a
$G$-space $M$, to each $\xi$ in the Lie algebra $\g$ of $G$, there
corresponds a fundamental vector field $\xi_M$ on $M$, given by
\begin{equation}
\label{eq:FVF}
	\xi_M(x) = \left.\frac{d}{dt}\exp t\xi \cdot x\right|_{t=0},
\end{equation}
for all $x \in M$.
This is a Lie algebra {\em anti}-homomorphism%
\footnote{Introducing a minus sign in the
definition (\ref{eq:FVF}) of the fundamental vector field would make
this a {\em bona fide} homomorphism of Lie algebras, but would not aid
in calculations or intuition.}  
$\g \to \LieAlg{X}(M)$; i.e., the
assignment $\xi \mapsto \xi_M$ is linear and satisfies, for $\xi, \eta
\in \g$,
\begin{equation}
\label{eq:FVF-anti}
	[\xi_M, \eta_M] = -[\xi, \eta]_M.
\end{equation}

An extremely useful construct for a $G$-manifold is the equivariant
de~Rham cohomology.  Note that the de~Rham complex $\Omega^*(M)$ of
smooth differential forms on $M$ has actions of $G$ and $\g$.  For $g
\in G$, we can pull-back by the diffeomorphism induced by $g$, and for
$\xi \in \g$, we have the Lie derivative and interior product
operators corresponding to the vector field $\xi_M$:
\begin{align*}
	\L_\xi = \L_{\xi_M} &\from \Omega^*(M) \longrightarrow
	\Omega^*(M); \\
	\i_\xi = \i(\xi_M) &\from \Omega^*(M) \longrightarrow
	\Omega^{*-1}(M); 
\end{align*}
The symmetric algebra $S(\g^*)$ has a natural action of $G$ given by
the coadjoint action on $\g^*$.  The {\em equivariant de~Rham
complex}\/ is given as a graded supercommutative algebra:
\[
	\Omega^k_G(M) \defeq \bigoplus_{2\ell + j = k}
	\left(\Omega^j(M) \tensor S^\ell(g^*)\right)^G.
\]
This complex can also be thought of as equivariant
$\Omega^*(M)$-valued polynomials on $\g$.  That is, $\alpha \in
\Omega^*_G(M)$ is a map $\g \to \Omega^*(M)$ which is polynomial in
any coordinates on $\g$ and satisfies
\[
	\alpha(\Ad_g\xi) = (g\inv)^*\alpha(\xi),
\]
for all $g \in G$ and $\xi \in \g$.  From this viewpoint, the Cartan
differential can be written 
\begin{equation}
\label{eq:d_G}
	\big(d_G\alpha\big)(\xi) = d\big(\alpha(\xi)) -
	\i_\xi\alpha(\xi).
\end{equation}
Of course, the cohomology of $\left(\Omega^*_G(M), d_G\right)$ is
called the equivariant de~Rham cohomology of $M$.  See \cite{GS-EC}
for a thorough treatment of equivariant de~Rham theory and how it
relates to classical Hamiltonian $G$-spaces.

\section{Definitions and (the) Examples}

\subsection{Moment Space and Moment Map}

Essentially, what we want a moment space to be is something that can
serve as the target of a moment map.  In the theory of Hamiltonian and
$q$-Hamiltonian $G$-spaces, we had notions of 
\begin{itemize}
\item
	The differential equation characterising a moment map;
\item
	Reduction of $G$-spaces;
\item
	Fusion; i.e., a multiplicative structure on the particular
	class of $G$-spaces;
\end{itemize}
and so on.  The first of these is the one we intend to axiomatize.

\begin{definition}
Let $G$ be a compact Lie group.
A {\em premoment space}\/ for $G$ is a pair $(P,\Xitwiddle)$, where $P$
is a $G$-manifold and $\Xitwiddle$ is a ``natural'' closed $G$-equivariant
three-form.  
\end{definition}

By ``natural'' we mean that $\Xitwiddle$ is functorial with respect to
inclusions.  Since 
\[
	\Omega^3_G(M) = \Omega^3(M)^G \oplus \left(\Omega^1(M) \tensor
	\g^*\right)^G,
\]
we can write $\Xitwiddle$ as
\[
	\Xitwiddle = \Xi + \tau,
\]
where $\Xi \in \Omega^3(P)^G$ is the invariant piece and $\tau \from \g
\to \Omega^1(P)$ is the equivariant piece.  It is a linear map.  The
condition that $\Xitwiddle$ be closed can be written as three equations:
\begin{align}
	d\Xi &= 0,		\\
\intertext{and for all $\xi \in \g$,}
	\i_\xi\Xi &= d\tau(\xi)		\\
\label{eq:Xitw-closed-3}
	0 &= \i_\xi\tau(\xi)		
\end{align}

This equivariant form allows us to define a moment map.
\begin{definition}
Let $M$ be a $G$-space and $P$ a premoment space for $G$.  $M$ is called
{\em $P$-momental}\/ if there exists an invariant two-form $\omega \in
\Omega^2(M)^G$ and equivariant map $\Phi \from M \to P$ such that
\begin{equation}
\label{eq:moment}
	d_G\omega = - \Phi^*\Xitwiddle.
\end{equation}
$M$ is called {\em $P$-Hamiltonian}\/ if in addition, for all $x \in M$,
\begin{equation}
\label{eq:min-deg}
	\ker\omega_x = 
	\setst{\xi_M(x)}{\xi \in \ker \tau_{\Phi(x)} \from
	\g \to T^*_{\Phi(x)}P}.
\end{equation}
\end{definition}

\begin{remarks}\ 
\begin{enumerate}
\item	
	The condition (\ref{eq:moment}) can be written less obliquely as
\begin{align}
\label{eq:mom1}
	d\omega &= - \Phi^*\Xi; \\
\label{eq:mom2}
	\iota(\xi_M)\omega &= \Phi^*\tau(\xi),
\end{align}
for all $\xi \in \g$.
\item
	For $p \in P$, $\tau_p$ is defined to be the linear map $\g \to
	T^*_pP$ which takes $\xi \in \g$ to the evaluation of the
	one-form $\tau(\xi)$ at the point $p$.  In light of
	(\ref{eq:mom2}), we have that for $p \in P$, the fundamental
	vector fields of all vectors in the kernel of $\tau_p$ must
	annihilate $\omega$.  Thus (\ref{eq:min-deg}) is a condition
	of minimal degeneracy.
\end{enumerate}
\end{remarks}

For a given $P$, the most immediate candidates for $P$-momental
$G$-spaces are the orbits of $G$.  These have a natural inclusion map
$i \from \O \to P$.  Indeed $\Xitwiddle$ induces a two-form on each
orbit $\O$.  If $p \in \O$, $T_p\O$ is spanned by $\setst{\xi_M(p)}{\xi
\in \g}$, and we define
\begin{equation}
\label{eq:omega-O}
	\omega_\O\left(\xi_P(p),\eta_P(p)\right) = 
	\tau(\xi)_p\left(\eta_P(p)\right).
\end{equation}
By (\ref{eq:Xitw-closed-3}), this form is well-defined and
alternating, and we immediately see that it satisfies (\ref{eq:mom2}).
$\omega_\O$ is characterised by the property that
\[
	\i_\xi\omega_\O = i^*\tau(\xi).
\]
We claim $\omega_\O$ is $G$-invariant, and this is a consequence of
the equivariance of $\tau$.  For, given $\xi, \eta \in \g$, and $g \in
G$, 
\begin{align*}
\left.g^*\omega_\O\right|_p\big(\xi_P(p),\eta_P(p)\big) 
	&= \left.\omega_\O\right|_{gp} \big(g_* \cdot \xi_P(p), 
		g_* \cdot \eta_P(p)\big) \\
	&= \left.\omega_\O\right|_{gp}\big((\Ad_g\xi)_P(gp), 
		g_* \cdot \eta_P(p)\big) \\
	&= \tau(\Ad_g\xi)_{gp}\big(g_* \cdot \eta_P(p)\big) \\
	&= \big((g\inv)^*\tau(\xi)\big)_{gp}
		\big(g_* \cdot \eta_P(p)\big) \\
	&= \tau(\xi)_p\left((g\inv)_*\cdot g_* \cdot \eta_P(p)\right) \\
	&= \left.\omega_\O\right|_p\big(\xi_P(p), \eta_P(p)\big).
\end{align*}
From the Cartan ``Magic Formula,''
\[
	0 = \L_\xi\omega_\O = d \i_\xi \omega_\O + i_\xi d\omega_\O
\]
we must have that
\begin{align*}
\i_\xi d\omega_\O 
	&= - d\i_\xi\omega_\O \\
	&= -di^*\tau(\xi) \\
	&= i^*d\tau(\xi) \\
	&= -i^* \i_\xi \Xi,
\end{align*}
which verifies the moment condition (\ref{eq:mom1}).  We have proved
the following:
\begin{proposition}
Let $P$ be a premoment space for $G$. Consider a $G$-orbit $\O \subset
P$ with two-form $\omega_\O$ given by (\ref{eq:omega-O}) and inclusion
map $i$.  Then $(\O, \omega_\O, i\from \O\to P)$ is a $P$-momental
$G$-space.
\end{proposition}
\begin{definition}
A premoment space $P$ is called a {\em moment space}\/ for $G$ if all
orbits $\O \subset P$ are $P$-Hamiltonian $G$-spaces with two-form given
by (\ref{eq:omega-O}) and moment map given by inclusion.
\end{definition}

This definition encompasses the heretofore known examples as we shall
now see.

\subsection{Examples of Moment Spaces}

The oldest and most well-known example of moment space is the dual $\g^*$
to the Lie algebra $\g$ of $G$.  Recall \cite{GS-ST} that a $G$-space $M$
is called {\em Hamiltonian}\/ (which we for the purposes of this
monograph will call {\em $\g^*$-Hamiltonian}) if there exists $\omega \in
\Omega^2(M)^G$ and a map $\Phi \from M \to \g^*$ equivariant with
respect to the coadjoint action of $G$ on $\g^*$ such that
\begin{align}
	d\omega &= 0 ; \\
	\iota(\xi_M)\omega &= d\tau(\xi);\\
	\ker \omega_x &= 0
\end{align}
for all $\xi \in \g$ and $x \in M$.  The form $\omega$ is otherwise
known as symplectic.  To see the equivariant three-form, we let
$\alpha \from \g \to C^\infty(\g^*)$ be the evaluation map
$\alpha(\xi)(\ell) = \left< \ell, \xi \right>$.  This is easily seen
to be an equivariant map, so $\alpha$ is an equivariant two-form.
$\Xitwiddle = \tau $ is just $d_G \alpha$, and therefore completely
equivariant (having no invariant part); for $\ell \in \g^*$ and
$\lambda \in T_\ell \g^*$,
\[
	\tau(\xi)_\ell(\lambda) = \left< \lambda, \xi \right>
\]
Each $\tau_\ell$ is obviously an isomorphism.  Let us denote by
$\theta_{\g^*}$ the $\g^*$-valued one-form which identifies each
$T_\ell \g^*$ with $\g^*$.  Then we can write
\[
	\tau(\xi) = \left< \theta_{\g^*}, \xi \right>.
\]

It is also well-known that coadjoint orbits $\O \subset \g^*$ are
symplectic manifolds and indeed Hamiltonian $G$-spaces.

The second example of moment space is the Lie group itself, as
examined in detail by Alekseev, Meinrenken, and Malkin.  Let $\g$ have
an invariant inner product $B$ attached to it (so $\G$ may be for
example the category of compact reductive groups).  Let $\theta =
\theta_G$ (respectively, $\thetabar = \thetabar_G$) in $\Omega^1(G,
\g)$ be the left- (respectively, right-)-invariant Maurer-Cartan form
on $G$.  So for $g \in G$,
\begin{align*}
	\theta_g &= L_{g\inv*} \from T_g G \to T_eG = \g; \\
	\thetabar_g &= R_{g\inv*} \from T_g G \to T_eG = \g.
\end{align*}
In a faithful matrix representation of $G$ (and this is the arena in
which it is most convenient to do calculations), $\theta = g\inv dg$ and
$\thetabar = dg\,g\inv$.  In this case, there is a closed bi-invariant
three-form which is canonical with respect to the pairing $B$:
\[
	\Xi = \frac{1}{12}B\big(\theta, [\theta, \theta]\big)
	= \frac{1}{12}B\big(\thetabar, [\thetabar, \thetabar]\big),
\]
where by the multiple appearances of $\theta$ or $\thetabar$ we mean the
skew-symmetrization of the corresponding multilinear map (this is a
convention we shall follow throughout).  $\Xi$ is closed due to the
invariance of $B$ and the {\em Cartan Structure Equations}
\begin{align}
	d\theta &= -\frac12[\theta,\theta]; \\
	d\thetabar &= -\frac12[\thetabar,\thetabar]
\end{align}

Let $G$ act on itself by conjugation.  Then for $\xi \in \g$, the
fundamental vector field on $G$ corresponding to $\xi$ is
\[
	\xi_G = \xi_L - \xi_R,
\]
where $\xi_L$ is the right-invariant vector field on $G$ generated by
the left translation action and {\em vice versa}\/ for $\xi_R$.  Then
\[
	\i(\xi_G)\Xi =
	\frac{1}{12}\i(\xi_L)B(\thetabar,[\thetabar,\thetabar]) 
	-\frac{1}{12}\i(\xi_R)B(\theta,[\theta,\theta]). 
\]
By invariance of $B$, this is
\[
	\frac14B(\xi,[\thetabar,\thetabar])
	-\frac14B(\xi,[\theta,\theta])
	= \frac12dB(\xi,\theta + \thetabar).
\]
So we may write $\tau(\xi) = \frac12B(\theta + \thetabar, \xi)$ as the
equivariant (with respect to the conjugation action) piece of the
equivariantly closed three-form $\Xitwiddle = \Xi + \tau$ on $G$.
Note that since $\thetabar_g = \Ad_g\theta_g$, and since $B$
identifies $\g^*$ with $\g$, we have $\ker \tau_g = \ker (1 + \Ad_g)$.

Then following \cite{AMM}, we say a $G$-space $M$ is {\em
$q$-Hamiltonian}\/ (henceforth, {\em $G$-Ham\-il\-ton\-i\-an}\/) if there
exists $\omega \in \Omega^2(M)^G$ and $\Phi \from M \to G$ equivariant
with respect to the conjugation action of $G$ on itself, such that
\begin{align}
	d\omega &= -\Phi^*\Xi ; \\
	\iota(\xi_M)\omega &= \frac12 \Phi^*
			B(\theta + \thetabar, \xi); \\
	\ker \omega_x &= \setst{\xi_M(x)}{\xi \in \ker(1 + \Ad_{\Phi(x)})},
\end{align}
for all $\xi \in \g$ and $x \in M$.
It is in fact this formulation of the moment and Hamiltonian conditions
that motivates the general definition.  As in the case of $\g^*$, the
canonical transitive $G$-Hamiltonian $G$-spaces are conjugacy classes
$\mathcal C \subset G$.

A final example arises from complexification.  The pairing $B$ extends
to a complex-bilinear pairing on the Lie algebra $\g^\C = \g \tensor
\C$, which we will call $B^\C$.  If $G$ is simply connected, there
exists a unique simply connected complex Lie group $G^\C$ with Lie
algebra $\g^\C$ containing $G$ as a subgroup.  Exponentiating the
conjugation automorphism gives us a map $\overline\blank \from G^\C \to G^\C$,
which singles out $G$ as its fixed-point set.  Define
\[
	P^\C = \set{h \in G^\C : \bar h = h\inv}.
\]
$P^\C$ is preserved by the conjugation action of $G$.
We can restrict canonical forms on $G^\C$ to $P^\C$.  
Let $\theta^\C$ and $\thetabar^\C$ be the Maurer-Cartan forms on
$G^\C$.  Set $\theta_P = \left.\theta^\C\right|_{P^\C}$ and likewise
define $\thetabar_P$.  Set
\[
	\Xi_P = \frac{1}{12} \Im
	B^\C\left(\thetabar_P, \left[\thetabar_P,
	\thetabar_P\right]\right)
\]
A $G$-space $M$ is called a {\em $q$-Hamiltonian $G$-space with
$P^\C$-valued moment map}\/ (henceforth, {\em $P^\C$-Hamiltonian
$G$-space}\/) if there exists $\omega \in \Omega^2(M)^G$ and $\Phi \from
M \to P^\C$ equivariant with respect to the conjugation action 
of $G$ on $P$ such that
\begin{align}
	d\omega &= -\Phi^*\Xi_P ; \\
\label{eq:mom2complex}
		\iota(\xi_M)\omega &= \frac{1}{2\sqrt{-1}}
		\Phi^*B(\theta_P + \thetabar_P, \xi); \\
	\ker \omega_x &= 0 ,
\end{align}
for all $\xi \in \g$ and $x \in M$.  
It can be shown that the pairing $B(\theta_P + \thetabar_P, \xi)$ is
purely imaginary, so the equation (\ref{eq:mom2complex}) makes sense.
Then $\tau = \frac{1}{2\sqrt{-1}}B(\theta_P + \thetabar_P, \blank)$,
and when seeking the kernel of $\tau_p$ for $p \in P^\C$, we find that
all the eigenvalues of $\Ad_p$ are positive.  Hence $\tau_p$ is an
isomorphism.  

These spaces are introduced in \cite{Alekseev} but named in \cite{AMM}.
In both, it is shown that $P^\C$-Hamiltonian $G$-spaces are in
one-to-one correspondence with $\g^*$-Hamiltonian $G$-spaces, due mainly
to the fact that with respect to $B$ there is a canonical diffeomorphism
$\g^* \iso P^\C$.

A different perspective on this example is through the category of
Poisson-Lie $G$-spaces, e.g., \cite{FL}.  In that category, the moment
map has as its target the {\em dual group}\/ $G^*$ of $G$ lying in
$G^\C$.  However, $\omega$ is no longer $G$-invariant, but instead
the $G$-action induces a Poisson map $G \times M \to M$.  

A crucial observation of these three moment spaces is that they appear
as coset spaces of a larger Lie group containing $G$ as a subgroup.
Indeed, for $G$ satisfying all reductivity and connectivity conditions
necessary: 

\begin{enumerate}
\item
Let $H_0 = G \rtimes \g^*$, so 
\[
	(g_1, \ell_1)\cdot(g_2, \ell_2) = (g_1 g_2, \ell_1 +
	\Ad^*_{g_1}\ell_2). 
\]
This is the right trivialization of the cotangent group $T^*G$.  The map
$j_0 \from H_0 \to \g^*$, $(g, \ell) \mapsto \ell$ gives an equivariant
diffeomorphism between the right coset space\footnote{%
If $G \subset H$ are Lie groups, we follow the convention that the {\em
right}\/ coset space refers to the space of orbits for the right action of
$G$ on $H$.  That is, $H/G$ is the space of all cosets $hG$ as $h$
ranges through $H$.}%
$H_0/G$ with left
$G$-action and $\g^*$ with the coadjoint $G$-action.

\item
Let $H_+ = G \times G$, with $G$ embedded as the diagonal subgroup
$\Delta(G)$.  The map $j_+\from G \times G \to G$, $(g_1, g_2) \mapsto g_1
g_2\inv$ gives an equivariant diffeomorphism of $H_+/\Delta(G)$ with $G$
equipped with its conjugation action.

\item
Let $H_- = G^\C$, which we have already defined.  The map $G^\C \to
G^\C$, $h \mapsto h \bar h\inv $ is seen to land in $P^\C$ and descends
to an equivariant diffeomorphism between $G^* = G^\C/G$ with its left
(``dressing'') action and $P^\C$ with the conjugation action. 
\end{enumerate}

Indeed, when we take the moment forms on $\g^*$, $G$, and $P^\C$ and
pull them back to the respective larger Lie groups, we find a way in
which they can all be expressed similarly.  This leads us to a general
recipe for providing moment spaces, which we will describe in the next
section. 

\section{Moment Spaces from Symmetric Pairs}

\subsection{Review of symmetric spaces, symmetric paris, and symmetric
Lie algebras}

A symmetric space is essentially a space with reflections through each
point.  Given such a Riemannian manifold $M$, we can consider the Lie
group $H$ (more precisely, $H$ is the connected component of this
group) of isometries of $M$, and the closed subgroup $G$ of all
isometries of $M$ fixing a given point.  Then $M \iso H/G$ are
$G$-equivariantly diffeomorphic.  This pair of groups is easier to
work with than $M$ alone, so we axiomatize the data

\begin{definition}[\cite{Helgason}]
Let $H$ be a Lie group and $G$ a closed subgroup.  $(H,G)$ is a
{\em symmetric pair}\/ if there exists an involutive automorphism 
$\sigma$ of $H$ (that is, $\sigma^2 = \Id$ while $\sigma\neq\Id$) such
that 
\begin{equation}
\label{eq:sigma-fixed}
	(H^\sigma)_0 \subseteq G \subseteq H^\sigma. 
\end{equation}
Here $H^\sigma$ is the subgroup of $H$ fixed by $\sigma$ and
$(H^\sigma)_0$ its identity component.

A symmetric pair $(H,G)$ is called {\em Riemannian}\/ if the image of
$\Ad_G$ in $\LieAlg{gl}(\h)$ is compact.
\end{definition}

\begin{remarks}
This definition is quite general, and could use some illumination.
\begin{enumerate}
\item
	Since in our context we are fixing $G$ and considering the
	different Lie groups $H$ such that $(H,G)$ is a symmetric
	pair, we can restrict our attention to connect groups $G$.  In
	this case (\ref{eq:sigma-fixed}) can be replaced by the more
	transparent identity $G = H^\sigma$.
\item
	The symmetric pair $(H,G)$ is Riemannian if in particular $G$
	is compact.  
\end{enumerate}
\end{remarks}

Let $\h$ be the Lie algebra of $H$ and $s = d\sigma_e$.  Then $(\h,s)$
has many linear properties analogous to $(H,G)$.  It is in fact an
example of the following object.

\begin{definition}
Let $\h$ be a Lie algebra over $\R$ and $s$ an involutive automorphism
of $\h$.  The pair $(\h,s)$ is called a {\em symmetric Lie algebra}.
$(\h,s)$ is called {\em orthogonal}\/ if $\g = \h^\sigma$ (the fixed
point subalgebra of $\h$) is compactly embedded in $\h$, and is called
{\em effective}\/ if $\g \intersect \LieAlg{z} = 0$, where $\LieAlg{z}$
is the center of $\h$.
\end{definition}

In particular, if $(H,G)$ is a Riemannian symmetric pair, the
associated symmetric Lie algebra $(\h,s)$ is orthogonal.  If $G$ is
semisimple, then $(\h,s)$ is effective.

\begin{lemma}
Let $(\h,s)$ be a symmetric Lie algebra.  Write
\[
	\h = \g \oplus \p
\]
for the decomposition of $\h$ into the subalgebra $\g$ fixed by
$\sigma$ and the -1 eigenspace $\p$ of $\sigma$.  Then
\[
\begin{array}{ccc}
	[\g, \g] \subset \g; &
	[\g, \p] \subset \p; &
	[\p, \p] \subset \g.
\end{array}
\]
\end{lemma}
\begin{proof}
$s$ is a Lie algebra homorphism.  Thus
\[
	s[\zeta_1,\zeta_2] = [s\zeta_1,s\zeta_2] 
\]
for all $\zeta_1$, $\zeta_2 \in \h$.
\end{proof}

Let $G$ be fixed and let $H_0$, $H_+$ and $H_-$ be the three Lie
groups given in Examples~1--3 above.  We will demonstrate that they
form symmetric pairs.

\begin{enumerate}
\item
	On $H_0 = G \rtimes \g^*$, the involution $\sigma_0$ is the map
	$(g,\ell) \mapsto (g, -\ell)$.  The corresponding symmetric Lie
	algebra is $\h_0 = \g \rtimes \g^*$ with Lie bracket and
	involution
\begin{align*}
	[\xi \rtimes \lambda, \eta \rtimes \mu]
	&=[\xi,\eta] \rtimes \left(\ad^*_\xi\mu -
	\ad^*_\eta\lambda\right) \\
	s_0(\xi \rtimes \lambda) &= \xi \rtimes (-\lambda).
\end{align*}
	The symmetric Lie algebra $\h_0$ is said to be of the
	{\em Euclidean type}.%
\footnote{We are using the same symbols as Helgason to label these
types, but the plus and minus signs are reversed.  Under our choice of
convention, $+$, $-$, and $0$ refer not only to the signs of the
sectional curvatures of the corresponding symmetric spaces but to the
signs of the eigenvalues in a certain linear transformation of $\g$
which will appear in the next section.  Ultimately, the permutation is
due to the fact that in the semisimple case, we will take the {\em
negative}\/ of the Killing form on $\g$ so as to get a {\em
positive}-definite inner product on $\g$.}
\item
	On $H_+ = G \times G$, $G$ is embedded as the diagonal
	$\Delta(G)$.  This subgroup is fixed by the involution
	$\sigma_+(g_1, g_2) = (g_2, g_1)$.  The corresponding
	symmetric Lie algebra is of course $\h_+ = \g \times \g$ with
	involution 
\[
	s_+(\xi_1 \times \xi_2) = \xi_2 \times \xi_1,
\]
	fixing the diagonal subalgebra $\LieAlg d(\g)$.  This
	symmetric Lie algebra is of the {\em compact type}.
\item
	On $H_- = G^\C$, we already know the conjugation 
	automorphism of $\g^\C$ and its exponentiated version.
	These single out the real forms of $\g^\C$ and $G^\C$,
	respectively.  
\[
	s_-(\xi + \sqrt{-1}\eta) = \xi - \sqrt{-1}\eta.
\]
	This symmetric Lie algebra is of the {\em noncompact type}.
\end{enumerate}
\subsection{Legendrian structure}

The symmetric pairs given in the last subsection are distinguished;
they have a pairing on the associated symmetric Lie algebra which
restricts to a nontrivial pairing between the different eigenspaces of
$s$.  We make the notion of this type of pairing precise with a
definition.  

\begin{definition}
Let $\LieAlg w$ be a vector space over a field $k$ with involution $s
\from \LieAlg w \to \LieAlg w$.  A {\em Legendrian form}\/ on
$\LieAlg w$ is a nondegenerate symmetric bilinear form $\Lambda \from
\LieAlg w \times \LieAlg w \to k$ with respect to which $s$ is
skew-symmetric, i.e., for all $w_1, w_2 \in \LieAlg w$,
\begin{equation}
\label{eq:L-s}
	\Lambda(sw_1, w_2) = - \Lambda(w_1, sw_2). 
\end{equation}
\end{definition}

Note that since $s^2 = \Id$, we can also write (\ref{eq:L-s}) as 
$\Lambda(sw_1, sw_2) = - \Lambda(w_1, w_2)$.

\begin{definition}
Let $(\h, s)$ be a symmetric Lie Algebra.  $\h$ is {\em Legendrian}\/
if it admits a Legendrian form which is invariant with respect to the
adjoint action of $\h$ on itself: for all $\zeta_1, \zeta_2, \zeta_3
\in \h$,
\begin{align}
\label{eq:L-ad}
	\Lambda(\ad_{\zeta_1}\zeta_2, \zeta_3) =
	-\Lambda(\zeta_2, \ad_{\zeta_1}\zeta_3),
\intertext{or, using symmetry,}
	\Lambda([\zeta_1, \zeta_2], \zeta_3) =
	\Lambda(\zeta_1, [\zeta_2,\zeta_3]).
\end{align}
\end{definition}

\begin{definition}
Let $(H,G)$ be a symmetric pair.  $H$ is {\em Legendrian} if the
associated symmetric Lie algebra $\h$ admits a Legendrian form which
is invariant with respect to the adjoint action of $H$ on $\h$.  That
is, in addition to (\ref{eq:L-s}) and (\ref{eq:L-ad}), we must have,
for all $\zeta_1, \zeta_2 \in \h$ and $g \in G$,
\begin{equation}
\label{eq:L-Ad}
	\Lambda(\Ad_h\zeta_1, \zeta_2) = 
	\Lambda(\zeta_1, \Ad_{h\inv}\zeta_2).
\end{equation}
\end{definition}

Let $\LieAlg w$ be a vector space with involution $s$, and write
$\LieAlg w = \LieAlg v \oplus \p$ as its decomposition into the +1 and
-1 eigenspaces of $s$.  Then a Legendrian form $\Lambda$ on $\LieAlg
w$ is zero when restricted to $\LieAlg v \times \LieAlg v$ and $\p
\times \p$.  That is, these eigenspaces are {\em isotropic} with
respect to $\Lambda$.  Since $\LieAlg v$ and $\p$ are algebraically
complementary, $\Lambda$ is a nondegenerate pairing between them, and
in particular $\LieAlg w$ is even-dimensional, isomorphic to $\LieAlg
v \oplus \LieAlg v^*$.  If in addition we are given an inner product
on $\LieAlg v$, we get a linear map $\LieAlg v \to \p$, which is
related to the Legendre transform.  Hence the attribution.

\begin{remark}
This also seems to be related to the concept of a Drinfeld algebra,
which is a Lie algebra structure on $\g \oplus \g^*$.  However, in
that formulation, $\g^*$ is required to be a subalgebra, whereas we
will always have $[\g^*,\g^*] \subset \g$.  Only in the least
interesting, Euclidean case do these definitions overlap.
\end{remark}

\begin{proposition}
Let $\g$ be a Lie algebra.
\begin{enumerate}[\upshape (a)]
\item
$\h_0 = \g \rtimes \g^*$ has Legendrian structure given by
\begin{equation}
	\Lambda_0(\xi_1 \rtimes \lambda_1, \xi_2 \rtimes \lambda_2) 
	= \left<\xi_1, \lambda_2\right> +\left<\xi_2,
	\lambda_1\right>. 
\end{equation}
$(H_0 = G \rtimes \g^*, G)$ is a Legendrian symmetric pair.
\item
Let $\g$ have an invariant inner product $B$.  $\h_+ = \g \times \g^*$ has
Legendrian structure given by
\begin{equation}
	\Lambda_+(\xi_1 \times \eta_1, \xi_2 \times \eta_2) 
	= \frac12\left(B(\xi_1, \xi_2) - B(\eta_1, \eta_2)\right)
\end{equation}
If $G$ is connected, $(H_+ = G \times G, \Delta(G))$ is a Legendrian
symmetric pair.
\item
Again assume $\g$ has an inner product $B$.  Then $B$ extends to a
$\C$-bilinear inner product on $\h_- = \g \tensor \C$.  $\h_-$ has a
Legendrian structure given by
\begin{equation}
	\Lambda_-(\zeta_1, \zeta_2)
	= \Im B^\C(\zeta_1, \zeta_2).
\end{equation}
If $G$ is simply connected, then $(H_- = G^\C, G)$ is a Legendrian
symmetric space.
\end{enumerate}
\end{proposition}
\begin{proof}
Clear.
\end{proof}

\subsection{The equivariant form}

The {\em Ansatz}\/ here is to show that given a Legendrian symmetric
space $(H,G)$, we can construct a moment space for $G$.  This space
will in fact be the quotient of $H$ by $G$.  We proceed in a series of
steps.

For the remainder of this section $(H,G)$ will be a Legendrian
symmetric pair with involution $\sigma$.  We write $\h = \g \oplus \p$
as the decomposition of the Lie algebra of $H$ into the eigenspaces of
$s = d\sigma|_e$.  The Legendrian structure on $\h$ will be denoted by
$\Lambda$.

Let $\theta$ be the left-invariant Maurer-Cartan form on $\h$.  Given
the involution $s$, we can decompose $\theta$ into its ``$\g$-part''
and ``$\p$-part''.  That is, write
\begin{align*}
	\gamma &= \frac{1+s}2 \theta; &
	\pi &= \frac{1-s}2 \theta;
\end{align*}
so $\gamma \in \Omega^1(H,\g)$ and $\pi \in \Omega^1(H,\p)$.  Let $j
\from H \to P$ be the quotient map.

\begin{proposition}
\label{prop:tau}
Define for $\xi \in \g$ a one-form
\begin{equation}
\label{eq:beta}
	\beta(\xi)_h = \Lambda(\xi, \Ad_h\pi) \in \Omega^1(H).
\end{equation}
Then
\begin{enumerate}[\upshape (a)]
\item
$\beta(\xi)$ is basic with respect to the right action of $G$ on $H$,
so there is a unique one-form $\tau(\xi) \in \Omega^1(P)$ such that
$j^*\tau(\xi) = \beta(\xi)$.
\item
The map $\xi \mapsto \beta(\xi)$ is equivariant with respect to the
left action of $G$ on $H$, so $\tau$ is an equivariant three-form on
$P$.
\item We have, for all $\xi \in \g$,
\begin{equation*}
	\i_\xi \tau(\xi) = 0.
\end{equation*}
\end{enumerate}
\end{proposition}
\begin{proof}
For $h \in H$, let $R_h$ and $L_h$ denote left and right
multiplication by $h$ as diffeomorphisms of $H$.  Since $R_g^*\theta =
\Ad_{g\inv}\theta$ and $\sigma(g) = g$, it follows that $R_g^*\pi =
\Ad_{g\inv}\pi$.  Then
\begin{align*}
	\big(R_g^*\beta(\xi)\big)_h 
	&= \Lambda(\xi, R_g^*\Ad_h\pi) \\
	&= \Lambda(\xi, \Ad_{hg}\Ad_{g\inv}\pi) \\
	&= \Lambda(\xi, \Ad_h\pi) = \beta(\xi)_h,
\end{align*}  
so $\beta(\xi)$ is right-invariant.  Morever, if $\eta_R(h) = L_{h*}
\eta$ is the fundamental vector field associated to the right action
corresponding to $\eta$, then $\theta(\eta_R) = \eta$.  Hence
$\pi(\eta_R) = 0$ and 
\begin{equation*}
	\beta(\xi)_h(\eta_R) = 0.
\end{equation*}
$\beta(\xi)$ is also right-horizontal, hence right-basic.  This proves
the first claim of the proposition.

For the second, note that $\theta$ and hence $\pi$ are left
$H$-invariant, so
\begin{align*}
	\big(L_{g\inv}^*\beta(\xi)\big)_h
	&= \Lambda(\xi, Ad_{g\inv h}\pi) \\
	&= \Lambda(\Ad_g\xi, \Ad_h\pi) = \beta(\Ad_g\xi)_h. 
\end{align*}

Finally, to prove the third claim, we will show that for $\xi \in \g$,
\begin{equation}
\label{eq:beta(xiL)=0}
	\i(\xi_L)\beta(\xi) = 0.
\end{equation}  
Indeed,
\begin{align}
\notag
	\beta(\xi)_h(\xi_L) 
	&= \Lambda\left(\xi, 
		\Ad_h\frac{\Ad_{h\inv} - \Ad_{\sigma(h\inv)}}2\xi
		\right) \\
\label{eq:AdH-invce}
	&= \frac12 \Lambda(\xi, \xi) -
	   \frac12 \Lambda(\Ad_{h\inv}\xi, \Ad_{\sigma(h\inv)}\xi) \\ 
	&= 0.  
\end{align}
Here the last step is from the $s$-skewness of $\Lambda$.  The
proposition is proved.
\end{proof}

\begin{remark}
It is only in (\ref{eq:AdH-invce}) that we used the full
$\Ad_H$-invariance of the pairing $\Lambda$.  In fact, the first two
claims of Proposition~\ref{prop:tau} can be proven with only a pairing
between $\g$ and $\p$ which is $\Ad_G$-invariant (note $\Ad_G$
preserves the decomposition $\h = \g \oplus \p$).  There is a unique
extension of such a pairing to an $s$-skew pairing of the full Lie
algebra, and to force the associated one-form to obey
(\ref{eq:beta(xiL)=0}) is essentially the assumption that this
extension is completely invariant.
\end{remark}

\subsection{The invariant form}
To summarize the last subsection, we have found an equvariant three-form
$\tau \in \Omega^1(H/G)$ that has the hope of being part of a closed
form.  We now complete this process by exhibiting the invariant piece.

\begin{proposition}
\label{prop:Xi}
Define $\Upsilon \in \Omega^3(H)$ by
\begin{equation}
\label{eq:Upsilon}
	\Upsilon = \frac13 \Lambda(\pi, [\pi,\pi]).
\end{equation}
Then
\begin{enumerate}[\upshape (a)]
\item
$\Upsilon$ is right $G$-basic and left $G$-invariant.  Hence there
exists a unique $\Xi \in \Omega^3(P)^G$ such that $\Upsilon = j^*\Xi$.
\item
\begin{equation}
\label{eq:dXi=0}
	d\Xi = 0.
\end{equation}
\item  For $\xi \in \g$,
\begin{equation}
\label{i(xi)Xi=dtau(xi)}
	i_\xi \Xi = d\tau(\xi).
\end{equation}
\end{enumerate}
\end{proposition}

We postpone the proof for a quick lemma.

\begin{lemma}
If $\theta = \gamma + \pi$ is the decomposition of $\theta$ relative
to that of $\h$, then
\begin{align*}
	d\gamma &= -\frac12([\gamma,\gamma] + [\pi,\pi]); \\
	d\pi    &= -[\gamma,\pi].
\end{align*}
\end{lemma}
\begin{proof}
This is an immediate consequence of the bracket identities for a
symmetric Lie algebra and the Cartan structure equations.
\end{proof}

\begin{proof}[Proof of Proposition~\ref{prop:Xi}]\ \\
\begin{enumerate}[\upshape (a)]
\item
This is proved similarly to the analogous claim in
Proposition~\ref{prop:tau}.  
\item
We are helped symbolically by the summation conventions on $\theta$,
$\gamma$ and $\pi$.  Note that by the Jacobi Identity
\begin{equation*}
	\big[\pi,[\pi,\pi]\big] =
	\big[\theta,[\theta,\theta]\big] =
	\big[\gamma,[\gamma,\gamma]\big] = 0
\end{equation*}
Thus,
\begin{align*}
	d\Upsilon 
	&= \frac13 d\Lambda(\pi, [\pi,\pi]) \\
	&= \Lambda(d\pi,[\pi,\pi]) \\
	&= -\Lambda([\gamma,\pi],[\pi,\pi]) \\
	&= -\Lambda\left(\gamma,\big[\pi,[\pi,\pi]\big]\right) = 0,
\end{align*}
so (\ref{eq:dXi=0}) is proved.
\item
Let $\xi \in \g$.  Then
\begin{align*}
	\i(\xi_L)\Upsilon 
	&= \frac13 \i(\xi_L)\Lambda(\pi,[\pi,\pi]) \\
	&= \Lambda\left(\pi(\xi_L),[\pi,\pi]\right) \\
	&= \Lambda(Ad_{h\inv}\xi,[\pi,\pi]).
\end{align*}
On the other hand,
\begin{align*}
	d\beta(\xi) &= d\Lambda(\xi, \Ad_h\pi) \\
	&= \Lambda(\xi, \Ad_h\ad_\theta\pi) 
		-\Lambda(\xi,\Ad_h [\pi,\gamma]) \\
	&= \Lambda(\xi, \Ad_h[\gamma + \pi, \pi])
		-\Lambda(\xi,\Ad_h [\pi,\gamma]) \\
	&= \Lambda(\xi, \Ad_h[\pi,\pi]).
\end{align*}
Thus (\ref{i(xi)Xi=dtau(xi)}) is true as well.
\end{enumerate}
\end{proof}

As an immediate consequence, we have:
\begin{theorem}
If $(H,G)$ is a Legendrian symmetric pair, the equivariant three-form
$\Xitwiddle = \Xi + \tau$ is equivariantly closed, thus giving $P =
H/G$ the structure of a premoment space for $G$.
\qed
\end{theorem}

\subsection{Nondegeneracy}

Along with the equivariant condition (\ref{eq:moment}), which we have
just satisifed for an arbitrary Legendrian symmetric pair, there is
the nondegeneracy (actually, minimal degeneracy) condition
(\ref{eq:min-deg}).  This condition is nearly independent of the
results preceding, but becomes essential when actually working with
$P$-Hamiltonian spaces (note that any $G$-space with equivariant form
0 is a premoment space for $G$).  Here we will use the nondegeneracy
of the pairing to satisfy nondegeneracy of $\tau$.

\begin{proposition}
\label{prop:O-Ham}
Let $(H,G)$ be a Legendrian symmetric space, and $\O$ an orbit of
$G$ in $P = H/G$.  Then $\O$ with two-form given by (\ref{eq:omega-O})
and moment map $i \from O \to P$ satifies
\begin{equation}
\label{eq:i-min-deg}
	\ker\omega_p = \setst{\xi_P(p)}{\xi \in \ker \tau_{p} \from
	\g \to T^*_{p}P}.
\end{equation}
Hence $\O$ is a $P$-Hamiltonian $G$-space.
\end{proposition}

Again, the theorem follows immediately: 
\begin{theorem}
Let $(H,G)$ be a Legendrian symmetric pair.  Then $P = H/G$ is a
moment space for $G$.
\qed
\end{theorem}

\begin{proof}[Proof of Proposition~\ref{prop:O-Ham}]
What we are attempting to prove is 
\begin{equation}
\label{eq:O-min-deg}
	\i_\xi\omega = 0 \iff
\begin{cases}
	\xi_P = 0 	&\text{or} \\
	\xi \in \ker\tau_p
\end{cases}
\end{equation}
Suppose that $\xi\in\ker\tau_p$, where $p=hG$.  This means that
\[
	0=\Lambda\left(\Ad_{h\inv}\xi,
		\frac{1-s}2\Ad_{h\inv}\eta\right) 
	=-\frac12(\Ad_{h\inv}\xi,\Ad_{\sigma(h\inv)}\eta)
\]
for all $\eta \in \g$.  Since $\g = \g^\Lambda$, we have that 
\begin{equation}
\label{eq:Ad_kxi-in-g}
	\Ad_{\sigma(h)h\inv}\xi \in \g.
\end{equation}
Write $k = \sigma(h)h\inv$ and note that $\sigma(k) = k\inv$.  Then
by (\ref{eq:Ad_kxi-in-g}) we must have
\[
	\Ad_{k\inv} \xi = \sigma\left(\Ad_k\xi\right) = \Ad_k\xi
\]
and therefore $\xi = \Ad_{k^2}\xi$ or 
\[
	\xi \in \ker (1 - \Ad_{k^2}).
\]
Now we have a direct-sum decomposition
\[
	\ker(1 - \Ad_{k^2}) = \ker(1 - \Ad_k) \oplus \ker(1 + \Ad_k).
\]
If $\xi$ is in the first summand, we have 
$\Ad_{h\inv}\xi = \Ad_{\sigma(h\inv)}\xi \in \g$ 
and therefore $\xi_{P}(p) = 0$.  On the
other hand, if $\xi$ is in the second summand we have
\[
	\Ad_{h\inv}\xi = - \Ad_{\sigma(h\inv)}\xi \in \p 
\]
and thus $\beta(\xi)_h = 0$.  Therefore (\ref{eq:O-min-deg}) is true.
\end{proof}

\section{A Decomposition Theorem}

We have shown how Legendrian symmetric pairs can give moment spaces.
We now attempt to show the extent to which the known examples are the
only ones.

Of course, if $G=G_1 \times G_2$ is a direct product of Lie groups,
and $P_1$ and $P_2$ are moment spaces for $G_1$ and $G_2$,
respectively, then $P_1 \times P_2$ with equivariant form
$\Xitwiddle_1 + \Xitwiddle_2$ is a moment space for $G$.  Thus we have
a way of ``building up'' moment spaces.  It is natural to try to go
the other way---i.e., to decompose.

\begin{theorem}
\label{thm:h-decomp}
Let $(\h,s)$ be an effective, orthogonal symmetric Lie algebra with
Legendrian structure $\Lambda$.  Then there exists a unique canonical
(up to isometry) decomposition
\begin{align}
\label{eq:h-decomp}
	\h &= \h_0 \oplus \h_+ \oplus \h_-; 
	&&\text{(direct sum of ideals)} \\ 
\label{eq:g-decomp}
	\g &= \g_0 \oplus \g_+ \oplus \g_-;
	&&\text{(direct sum of ideals)} \\ 
\label{eq:p-decomp}
	\p &= \p_0 \oplus \p_+ \oplus \p_- 
	&&\text{(direct sum of subspaces)}
\end{align}
such that, with the induced symmetric and Legendrian structures given
by restriction
\begin{align*}
	\h_0 &= \g_0 \oplus \p_0 \iso \g_0 \rtimes \g_0^*; \\
	\h_+ &= \g_+ \oplus \p_+ \iso \g_+ \times \g_+; \\
	\h_- &= \g_- \oplus \p_- \iso \g_- \tensor \C. \\
\end{align*}
\end{theorem}
\begin{remark}
This is very similar to the decompostion of effective, orthogonal
symmetric Lie algebras into a sum of Euclidean, compact, and
noncompact pieces given in \cite[Ch.~V, Theorem~1.1]{Helgason}.
\end{remark}

We will prove this in a sequence of lemmata.  To begin with, let $B$
be the negative of the Killing form on $\h$.  Then $B$ is
positive-definite on $\g$, $\ad_\h$-invariant, and also $s$-invariant:
\begin{align*}
	B(s\zeta_1,s\zeta_2) &= B(\zeta_1,\zeta_2), \\
\intertext{or,}
	B(s\zeta_1,\zeta_2) &= B(\zeta_1,s\zeta_2),
\end{align*}
for all $\zeta_1$, $\zeta_2 \in \h$.  Define $J \from \h \to \h$ by
\begin{equation}
\label{eq:J}
	\Lambda(J\zeta_1, \zeta_2) = B(\zeta_1, \zeta_2).
\end{equation}

\pagebreak[2]
\begin{lemma}\label{lem:J}\ \\
\begin{enumerate}[\upshape (a)]
\item \nopagebreak
\label{lem:J:a}
	$J$ commutes with the adjoint action of $\h$ on itself.  That
	is, for all $\zeta \in \h$,
\begin{align*}
	J \comp \ad_\zeta &= \ad_\zeta \comp J; \\
\intertext{or, for all $\zeta_1$ and $\zeta_2$,}
	J[\zeta_1,\zeta_2] &= [J\zeta_1, \zeta_2];
\end{align*}
\item
\label{lem:J:b}
	$J$ anti-commutes with $s$:  $J\comp s = - s \comp J$.  So $J$
	takes $\g$ into $\p$ and {\em vice-versa}.
\item
\label{lem:J:c}
	$J$ is self-adjoint with respect to $\Lambda$;.
\item
\label{lem:J:d}
	$J|_\g$ is a vector-space isomorphism $\g \iso \p$.
\end{enumerate}
\end{lemma}
\begin{proof}
Parts (\ref{lem:J:a})~and~(\ref{lem:J:b}) are straightforward.  Part
(\ref{lem:J:c}) is a simple consequence of symmetry of $B$ and
$\Lambda$.  Part (\ref{lem:J:d}) follows from the fact that $B$ is
positive definite on $\g$.
\end{proof}

It follows that $J^2$ is an endomorphism of $\g$.  By Lemma~\ref{lem:J},
Part (\ref{lem:J:c}), $J^2$ is self-adjoint.  Therefore, $\g$ has an
orthonormal (with respect to $B$) basis of eigenvectors with real
eigenvalues.  By rescaling $\Lambda$ by positive constants, we may
assume all nonzero eigenvalues are $\pm1$.  Put
\begin{align*}
	\g_0 &= \ker J^2;\\
	\g_+ &= \text{$+1$ eigenspace of $J^2;$} \\
	\g_- &= \text{$-1$ eigenspace of $J^2;$}
\end{align*}

\begin{lemma}
\label{lem:g-decomp}
$\g = \g_0 \oplus \g_+ \oplus \g_-$ is a direct sum decomposition of
ideals orthogonal with respect to $B$.  In fact,
\begin{equation}
\label{eq:g-comm}
	[\g_0, \g_+] =
	[\g_0, \g_-] =
	[\g_-, \g_+] = 0.
\end{equation}
\end{lemma}

\begin{proof}
Let $\xi$, $\eta \in \g$ have eigenvalues $\alpha$ and $\beta$,
respectively.  Then
\[
	\alpha B(\xi, \eta) = 
	B(J^2\xi, \eta) = B(\xi, J^2\eta) = \beta B(\xi,\eta).
\]
Hence $\alpha=\beta$ or $B(\xi,\eta) = 0$.  Also, 
\[
	\alpha[\xi, \eta] =
	[J^2\xi, \eta] =
	J^2[\xi,\eta] =
	[\xi, J^2\eta] = \beta[\xi, \eta],
\]
so the summands are ideals, and the commutation relations
(\ref{eq:g-comm}) hold.
\end{proof}

Now let
\begin{align*}
	\p_0 &= J\g_0	&\h_0 &= \g_0 \oplus \p_0; \\
	\p_+ &= J\g_+	&\h_+ &= \g_+ \oplus \p_+; \\
	\p_- &= J\g_-	&\h_- &= \g_- \oplus \p_-; 
\end{align*}
and note that (\ref{eq:p-decomp}) and (\ref{eq:h-decomp}) hold {\em
modulo} the proof that the summands of $\h$ are ideals.  We prove this
next.

\begin{lemma}
The following commutation relations hold:
\begin{equation}
\label{eq:p-comm}
	[\p_0, \p_+] =
	[\p_0, \p_-] =
	[\p_-, \p_+] = 0.
\end{equation}
Furthermore,
\[
\begin{array}{ccc}
	[\p_0, \p_0] = 0; &
	[\p_+, \p_+] \subset \g_+; &
	[\p_-, \p_-] \subset \g_-.
\end{array}
\]
\end{lemma}
\begin{proof}
Any element of $\p$ can be written as $J\xi$ for $\xi \in \g$.  Let
$\xi$ have $J^2$-eigenvalue $\alpha$ and $\eta$ have eigenvalue
$\beta$.  Then since
\[
	[J\xi,J\eta] = J^2[\xi,\eta] = \alpha[\xi,\eta] = \beta[\xi,\eta],
\]
we can apply Lemma~\ref{lem:g-decomp} and get the desired result.
\end{proof}

All we need to do now is supply the isomorphisms.  For $v \in
\p_0$, there is a unique covector $\i_\eta \Lambda \in \g_0^*$ (that
is, $\Lambda$ induces an isomorphism $\p_0 \iso \g_0^*$).  Define
\begin{align*}
	T_0 \from \h_0 &\longrightarrow \g_0 \rtimes \g_0^* \\
	(\xi, v) &\longmapsto \xi \rtimes \i_\eta \Lambda;
\end{align*}
$T_0$ is obviously bijective.  Since $[\p_0, \p_0] = 0$, it follows
that $T_0$ is a homomorphism of symmetric Lie algebras, and thus an
isomorphism.

Define another bijection
\begin{align*}
	T_+ \from \h_+ &\longrightarrow \g_+ \times \g_+ \\
	(\xi, J\eta) &\longmapsto \frac12(\xi + \eta, \xi-\eta);
\end{align*}
Given the fact that $[J\xi,J\eta] = [\xi,\eta]$ on this subalgebra,
$T_+$ is an isomorphism.  Finally, define
\begin{align*}
	T_- \from \h_- &\longrightarrow \g_- \tensor \C; \\
	(\xi, J\eta) &\longmapsto \xi + \sqrt{-1}\eta;
\end{align*}
Since this time $[J\xi,J\eta] = -[\xi,\eta]$, $T_-$ is also an
isomorphism.  This proves the theorem.

\begin{theorem}
Let $(H,G)$ be a Legendrian symmetric pair, 
with $G$ connected and semisimple, and $H$ simply connected.  
Then the moment space $P = H/G$ has
a decomposition
\[
	P = P_0 \times P_+ \times P_-
\]
and $G$ has a decomposition
\[
	G = G_0 \times G_+ \times G_-
\]
such that $P_0$ is a moment space for $G_0$ isomorphic to $\g_0^*$,
$P_+$ is a moment space for $G_+$ isomorphic to $G_+$, and $P_-$ is a
moment space for $G_-$ isomorphic to $P^\C(G_-)$.
\end{theorem}
\begin{proof}
Cf. \cite[Ch.~V, Prop.~4.2]{Helgason}
It follows from the homotopy exact sequence for the fibration
$G \to H \to H/G$ that if $H$ is simply connected and $G$ is
connected, then $H/G$ is simply connected.\footnote{%
Recall that for a topological group $G$, $\pi_0$ has a group structure.
In fact, $\pi_0(G) = G/G_0$, where $G_0$ is the connected component of
the identity.  Thus the exact sequence is one of groups even on the
$\pi_0$ level.}

Since $G$ is semisimple, $\h$ is an effective, orthogonal, Legendrian, 
symmetric Lie algebra  Therefore, we can compose $\h$ as in 
Theorem~\ref{thm:h-decomp} into $\h = \h_0 \oplus \h_+ \oplus
\h_-$.  Let $H_0 \times H_+ \times H_-$ be the
corresponding decomposition of $H$.  Likewise $\g$ decomposes and we
can write $G = G_0 \times G_+ \times G_-$.  Then
\[
	P = H/G 
	= \frac{H_0 \times H_+ \times H_-}{G_0 \times G_+ \times G_-} 
	= H_0/G_0 \times H_+/G_+ \times H_-/G_-.
\]
\end{proof}

\section{Conclusions}

Of course, we have said nothing about further axioms of a moment space
we might want to assume, in order to find the correct analogues of
fusion and reduction, etc.  It is quite possible that with the
addition of certain axioms, a moment space {\em must}\/ be a globally
symmetric space, and therefore must have an associated symmetric
pair.  This would complete the classification of moment spaces.

\bibliographystyle{habbrv} 
\nocite{AMM} 
\nocite{Helgason}
\nocite{Alekseev} 
\bibliography{math,moment}

\end{document}